\documentclass[12pt,a4paper]{amsart}
\usepackage{graphicx}
\usepackage{epstopdf}
\usepackage{amssymb}
\usepackage{amsthm}
\usepackage{amsmath}
\usepackage{pstcol, pst-node}
\usepackage{mathrsfs}
\usepackage{amscd}
\usepackage{enumerate}
\usepackage[all]{xy}

\newtheorem{thm}{Theorem}[section]

\newtheorem{prop}[thm]{Proposition}

\newtheorem{cor}[thm]{Corollary}

\newtheorem{lem}[thm]{Lemma}

\theoremstyle{definition}
\newtheorem{df}[thm]{Definition}
\newtheorem{definition}[thm]{Definition}
\newtheorem{rem}[thm]{Remark}

\newcommand{\dge}{\rotatebox[origin=c]{45}{$\ge$}}
\newcommand{\uge}{\rotatebox[origin=c]{315}{$\ge$}}

\newcommand{\omegatilde}{\widetilde{\omega}}

\newcommand{\bC}{\mathbb{C}}

\newcommand{\bP}{\mathbb{P}}

\newcommand{\frakX}{\mathfrak{X}}

\title[Tropical curves and Gromov-Witten invariants]{Toric degenerations, tropical curves and Gromov-Witten invariants
 of Fano manifolds}
\author{Takeo Nishinou}
\date{}

\begin{document}

\maketitle

\begin{abstract}
In this paper, we give a tropical method for computing Gromov-Witten 
 type invariants
 of Fano manifolds of special type.
This method applies to those Fano manifolds which admit toric degenerations
 to toric Fano varieties with singularities allowing small resolutions.
Examples include (generalized) flag manifolds of type A, and some moduli space
 of rank two bundles on a genus two curve.
\end{abstract}
\section{Introduction}
Since its appearance, 
 there are various kinds of applications of the ideas of
 tropical geometry to problems in classical geometry.
But in many cases the ambient spaces are toric varieties,
 and the number of applications
 of tropical geometry to problems in
 non-toric varieties are not large.
In this paper, we try to extend the applicability of
 tropical geometry in some direction, namely the enumerative problems.
More precisely, 
 in this paper we give a method to compute Gromov-Witten 
 type invariants of Fano manifolds of special
 type, in terms of counting of tropical curves.

Our calculation is based on toric degenerations of Fano manifolds.
Given such degenerations, we can compare holomorphic curves 
 between Fano manifolds and (singular) toric varieties.
On the other hand, the correspondence between such curves
 and tropical curves is shown in \cite{M,N2,NS}.
Thus, we can count curves in Fano manifolds by
 counting appropriate tropical curves.

The content of this paper is as follows.
After a short review of toric degenerations in Section 2, 
 we try to define Gromov-Witten type invariants via tropical method in Section 3.
In fact, we will define two types of invariants.
The first one (Theorem \ref{tropical}) is directly related to Gromov-Witten invariants
 and defined when the incidence conditions (Subsection \ref{subsec:incidence})
 satisfy suitable assumptions (Assumption 1 in Subsection \ref{subsec:incidence}
 and the transversality assumption in Lemma \ref{transverse}).
In particular, this invariant depends on the incidence conditions only through
 their homology classes.
The second one (Theorem \ref{tropical2}) is the analogue of the
 one considered in 
 \cite{NS}, and regarded as the 'relative part' of the Gromov-Witten invariants.
A priori it might not be a homological invariants, but it gives a lower bound to 
 the corresponding Gromov-Witten invariants
 (see Remark \ref{rem:GW}).
This second invariant is defined under only Assumption 1, and contains the first 
 invariant as a special case.

In the last section, we give two examples of 
 the calculation of these invariants for varieties of particular interest
 (these are both the first type invariants, so 
 we calculate actual Gromov-Witten invariants).
In these cases, in principle we can calculate all 
 genus zero Gromov-Witten invariants
 by tropical method.
In particular, we can calculate Gromov-Witten invariants including odd 
 cohomology classes.
Also, by the homological invariance of Gromov-Witten invariants, 
 an interesting combinatorial invariance of the counting number of tropical curves
 is observed, which cannot be deduced
 purely by tropical method (Remark \ref{rem:inv}).\\

\noindent
{\bf Acknowledgments.}
This research was inspired by the conversations with
 my collaborators Yuichi Nohara and Kazushi Ueda.
I would like to thank them for these conversations. 
I am supported by Grant-in-Aid for Young Scientists
(No. 19740034).

\section{Toric degeneration}
 \label{sc:toric_deg}
Let us first recall the definition of toric degenerations
 of projective manifolds.
\begin{df} \label{def:toric_degeneration}
Let $(X, \omega)$ be a projective manifold $X$
with a K\"{a}hler form $\omega$. 
A {\em toric degeneration} $(\mathfrak X, \widetilde \omega)$
 of $(X, \omega)$
 is a flat family
$
 f : \frakX \to B
$
of complex varieties over a connected complex variety $B$
 with two distinguished points $p_0, p_1$, and
a K\"{a}hler form $\omegatilde$ on $\frakX$
 (defined on the smooth locus), such that
$
 (X_1, \omega_1)
$
is isomorphic to $(X, \omega)$
as a K\"{a}hler manifold, and
$
 (X_0, \omega_0)
$
is a toric variety
with a torus invariant K\"{a}hler form. 
Here $(X_0, \omega_0)$, $(X_1, \omega_1)$ are the restrictions of
 $(\mathfrak X, \widetilde \omega)$ to the 
 fibers over $p_0, p_1$, respectively.

If $(X_0, \omega_0)$ is a toric Fano variety, then we call
 $(\mathfrak X, \widetilde \omega)$ a toric Fano degeneration
  of $(X, \omega)$. 
  
We note that, given such a family, if there is a holomorphic disk
 $\eta: D(2) = \{z\in \Bbb C\;|\; |z|<2\} \to B$ in the base space $B$, 
 with $\eta(0) = p_0$ and $\eta(1) = p_1$, there is a natural map, the 
 \emph{gradient Hamiltonian flow} (\cite{R}, see also \cite{NNU})
\[
\phi_{grH,t}: X_t\to X_0,\;\; t\in [0, 1]
\]
 which is a diffeomorphism away from the singular locus of $X_0$,
 and keeps the K\"ahler forms (also away from the singular points).
Here $X_t$ is the fiber over $\eta(t)$.
We write $\phi_{grH, 1}: X_1\to X_0$ by $\phi_{grH}$.
\end{df}
\begin{df}
Let $X$ be a (singular) toric variety defined by a
 fan $\Sigma$ in $\Bbb R^n$.
We say that $X$ allows a small resolution when 
 there is a refinement $\Sigma'$ 
 of $\Sigma$ without adding a new ray
 such that the toric variety associated to $\Sigma'$
 is nonsingular. 
\end{df}
An example of a class of Fano varieties having degenerations
 to toric Fano varieties allowing small resolutions
 is provided by flag manifolds of type A (including partial
 flag manifolds), see \cite{NNU}.
\begin{rem}\label{rem:rem}
In \cite{NNU}, we developed the notion of \emph{toric degeneration
 of integrable systems}.
Using this notion and the methods
in \cite{N1}, the results in this paper
 can be extended to the counting of holomorphic 
 disks with Lagrangian boundary conditions.
Note that flag manifolds of type A allow degenerations
 of integrable systems (\cite{NNU}).

Also, using the method of \cite{N2}, we can partially
 extend the results of this paper to
 higher genus tropical curves.
In fact, combining the methods of \cite{N1} and \cite{N2}, 
 we can even deal with curves with any genus and any number of boundary 
 components.
\end{rem}

\section{Computing Gromov-Witten type
 invariants of Fano manifolds by tropical method}
In this section, we explain that when we have a Fano manifold which has
 a degeneration to toric Fano variety \emph{allowing a small resolution}
 in the sense defined in
 Definition \ref{def:toric_degeneration},
 we can count the number of appropriate holomorphic curves
 in the Fano manifold.
This will be performed with the help of tropical method developed in
 the case of toric variety (\cite{M,N2,NS}).
Namely, in \cite{M,N2,NS}, an equality between suitable counts of tropical curves in 
 affine spaces and those of holomorphic curves in toric varieties are shown.
The latter is an analogue of Gromov-Witten invariants.

When there is a toric degeneration of a Fano manifold,
 we may compare the curves in the Fano manifold and 
 those in its degeneration, 
 that is, the toric variety.
When this comparison is effective, we may compute the number of 
 holomorphic curves in Fano manifold by counting appropriate tropical curves.
We will show this can be performed when the toric variety have
 a small resolution.

\subsection{Tropical curves}
First we recall some definitions about tropical curves, 
 see \cite{M,N2,NS}
 for more information.
 Let $\overline \Gamma$ be a weighted, connected finite graph.
Its sets of vertices and edges are denoted 
\[
\overline \Gamma^{[0]},
\overline \Gamma^{[1]}.\]
Then we write the weight function by 
\[w_{\overline \Gamma} : 
  \overline \Gamma^{[1]} \to \Bbb N \setminus \{ 0 \}.
  \]
An edge $E \in \overline \Gamma^{[1]}$ has adjacent vertices
 \[\partial E = \{ V_1, V_2 \}.\]
Let 
\[\overline \Gamma^{[0]}_{\infty} \subset \overline \Gamma^{[0]}\]
 be the set of one-valent vertices.
We set 
\[\Gamma = \overline \Gamma \setminus \overline\Gamma^{[0]}_{\infty}.\]
Non-compact edges of $\Gamma$ are called \emph{unbounded edges}.
Let $\Gamma^{[1]}_{\infty}$ be the set of unbounded edges.
Let 
\[
\Gamma^{[0]}, \Gamma^{[1]}, w_{\Gamma}
\]
 be the sets of vertices and edges of $\Gamma$ and the weight function
 of $\Gamma$ (induced from $w_{\overline\Gamma}$ in an obvious way),
 respectively.
Let $N$ be a free abelian group of rank $n\geq 2$
 and $N_{\Bbb R} = N\otimes_{\Bbb Z}\Bbb R$.

\begin{definition}\label{trop}
A \emph{parameterized trivalent tropical curve} in $N_{\Bbb R}$ is a proper map
 $h : \Gamma \to N_{\Bbb R}$ satisfying the following conditions.
\begin{enumerate}[(i)]
\item $\Gamma$ is a trivalent graph.
\item For every edge $E \subset \Gamma$ the restriction $h \big|_E$
 is an embedding with the image $h(E)$ 
 contained in an affine line with rational slope.
\item For every vertex $V \in \Gamma^{[0]}$, the following \emph{balancing
 condition} holds.
 Let $E_1, \dots, E_m \in \Gamma^{[1]}$ be the edges adjacent to $V$ and
 let $u_i \in N$ be the primitive integral vector emanating from $h(V)$
 in the direction of $h(E_i)$.
 Then
\begin{equation}
\sum_{j=1}^m w(E_j)u_j = 0.
\end{equation}
\end{enumerate}
\end{definition}
\emph{In this paper, we always assume that the graph $\Gamma$ is trivalent,
 so hereafter a (parameterized) tropical curve means a (parameterized)
 trivalent tropical curve.}

An isomorphism of parameterized tropical curves $h : \Gamma \to N_{\Bbb R}$ and 
 $h' : \Gamma' \to N_{\Bbb R}$ is a homeomorphism $\Phi : \Gamma \to \Gamma'$
 respecting the weights and $h = h' \circ \Phi$.
A \emph{tropical curve} is an isomorphism class of parameterized 
 tropical curves.
The \emph{genus} of a tropical curve is the first Betti number of $\Gamma$.
A \emph{rational} tropical curve is a tropical curve of genus zero.

The set of \emph{flags} of $\Gamma$ is 
\[
F\Gamma = \{(V, E) \big|
     V \in \partial E \}.
     \]
By (i) of the definition we have a map
\[
u : F\Gamma \to N
\] sending a flag $(V, E)$
 to the primitive integral vector $u_{(V, E)} \in N$
 emanating from $V$ in the direction of $h(E)$. 

An \emph{l-marked tropical curve} 
 is a tropical curve $h : \Gamma \to N_{\Bbb R}$
 together with a choice
 of $l$ not necessarily distinct
 edges 
\[
\bold{E} = (E_1, \dots, E_l) \in (\Gamma^{[1]})^l.
\]

The \emph{combinatorial type}
 of an $l$-marked tropical curve $(\Gamma, \bold{E}, h)$
 is the marked graph $(\Gamma, \bold{E})$ 
 together with the map $u : F\Gamma \to N$. 

The \emph{degree} of a type $(\Gamma, \bold{E}, u)$
 is a function $\Delta: N \setminus \{ 0 \}
  \to \Bbb N$
 with finite support defined by 
\begin{equation*}
 \Delta(\Gamma, u)(v):= \sharp \{ (V, E) \in F\Gamma |
    E \in \Gamma^{[1]}_{\infty}, w(E)u_{(V, E)} = v \} 
\end{equation*}
Let 
\[
e = |\Delta| = \sum_{v\in N\setminus\{0\}}\Delta(v).
\]
This is the same as the number of unbounded edges 
 of the graph $\Gamma$.
It is known that the space $\mathfrak T_{(\Gamma, \bold{E}, u)}$
 of marked tropical curves of given type
 $(\Gamma, \bold{E}, u)$ is, if non-empty, 
 a manifold with boundary (in fact, a convex polytope)
 of dimension larger than or equal to
\[
e+(n-3)(1-g)-ov(\Gamma).
\]
Here 
 $ov(\Gamma)$ is the \emph{overvalence} of $\Gamma$
 defined by
\[
ov(\Gamma) = 
\sum_{V\in \Gamma^{[0]}}(\sharp\{E\in\Gamma^{[1]}|(V, E)\in F\Gamma\}-3).
\]
It is also known that if $\Gamma$ is rational, the 
 equality 
\[
\dim \mathfrak T_{(\Gamma, \bold{E}, u)}=e+(n-3)(1-g)-ov(\Gamma)
\]
 holds.
\subsection{Incidence conditions}\label{subsec:incidence}
To fix the counting problem, we have to define {\it incidence conditions}
 for tropical curves and holomorphic curves.
We 
 recall some terminologies concerning incidence conditions from \cite{NS}.
See \cite{NS} for more details.
We formulate them for any genus in view of Remark \ref{rem:rem},
 although in this paper we almost always treat genus zero case.
\subsubsection{Incidence conditions for tropical curves}\label{subsubsec:trop}
We begin with the case of tropical curves.
\begin{definition}
For ${\bf d} = (d_1, \dots, d_l)\in \Bbb N^l$, an \emph{affine constraint} of codimension 
 ${\bf d}$ is an $l$-tuple ${\bf A} = (A_1, \dots, A_l)$
 of affine subspaces
 $A_i\subset N_{\Bbb R}$, defined over rational numbers,
 with 
\[
\dim A_i = n-d_i-1.
\]
An $l$-marked tropical curve $(\Gamma, {\bf E}, h)$ \emph{matches} the 
 affine constraint ${\bf A}$ if
\[
h(E_i)\cap A_i\neq \emptyset,\;\; i = 1, \dots, l.
\]
\end{definition}

Let us fix a degree $\Delta: N\setminus\{0\}\to\Bbb N$.
Now let 
\[
{\bf L} = (L_1, \dots, L_l)
\] be a set of linear subspaces of $N_{\Bbb Q}$,
 with ${\rm codim}\, L_i = d_i+1$.
Then the elements
\[
{\bf A} = (A_1, \dots, A_l), \;\; A_i\in N_{\Bbb Q}/L_i
\]
 define affine constraints.
\begin{definition}(Definition 2.4, \cite{NS})
Fix the genus $g$ and
 a degree $\Delta\in Map(N\setminus\{0\}, \Bbb N)$ and write $|\Delta| = e$
 as before.
An affine constraint ${\bf A} = (A_1, \dots, A_l)$ of codimension
 ${\bf d} = (d_1, \dots, d_l)$
  is \emph{general} for $\Delta$ and $g$ if
\[
\sum_{i =1}^l d_i = e+(n-3)(1-g),
\]
 and if any $l$-marked tropical curve $(\Gamma, {\bf E}, h)$ of genus g and degree
 $\Delta$ matching ${\bf A}$ satisfies the following:
\begin{enumerate}
\item $h(\Gamma^{[0]})\cap \bigcup_i A_i = \emptyset$.
\item $h$ is an embedding for $n>2$.
For $n = 2$, $h$ is injective on the subset of vertices, and for any 
 $x\in h(\Gamma)$, the inverse image is at most two points, and
 such $x\in h(\Gamma)$ with $\sharp\{h^{-1}(x)\} = 2$ is finite.
\end{enumerate}
\end{definition}
\begin{prop}(Proposition 2.4, \cite{NS})
Fix the genus $g$ and a degree $\Delta$ as above.
Let $\mathfrak A: = \prod_{i=1}^l N_{\Bbb Q}/L(A_i)$ be the space of 
 affine constraints of codimension ${\bf d} = (d_1, \dots, d_l)$
 such that $\sum_{i=1}^l d_i= e+(n-3)(1-g)$.
Then the subset
\[
\mathfrak Z:= \{{\bf A'}\in\mathfrak A\;|\; \text{${\bf A}'$ is non-general for
 $\Delta$ and $g$}\}
\]
 is nowhere dense in $\mathfrak A$.
\end{prop}
\proof
This is proved for $g = 0$ in \cite{NS}. 
The argument extends without little changes to any genus, and we omit it.\qed \\

For a marked tropical curve
 $(\Gamma, \bold{E}, h)$ matching the constraints $\bold{A}$,
 we have other important numbers, \emph{weight} and 
 \emph{index}.
The weight is defined by local data of the abstract graph
 $\Gamma$ (weights and markings of the edges)
 as 
\[
 w(\Gamma, \bold{E}) = \prod_{E\in \Gamma^{[1]}\setminus\Gamma^{[1]}_{\infty}}
  w_{\Gamma}(E)\cdot \prod_{i=1}^l w_{\Gamma}(E_i).
\]
There are two kinds of (lattice) indices, written as 
 $\mathfrak D(\Gamma, \bold{E}, h, \bold{A})$ and 
 $\delta_i(\Gamma, \bold{E}, h, \bold{A})$, 
 respectively (see Section 8, \cite{NS}).
The index $\mathfrak D(\Gamma, \bold{E}, h, \bold{A})$ is 
 defined as the index of the inclusion of the lattices (\cite{NS}, Proposition 5.7):
\[\begin{array}{rll}
 Map(\Gamma^{[0]}, N)&\to&\prod_{E\in \Gamma^{[1]}\setminus\Gamma_{\infty}^{[1]}}
 N/\Bbb Z u_{(\partial^- E, E)}\times\prod_{i=1}^l N/(\Bbb Qu_{(\partial^-
  E_i, E_i)}+L(A_i))\cap N,\\
h&\mapsto &\left((h(\partial^+E)-h(\partial^-E))_E, (h(\partial^-E_i))_i\right).
\end{array}
\] 
Here $\partial^{\pm}: \Gamma^{[1]}\setminus\Gamma^{[1]}_{\infty}\to\Gamma^{[0]}$
 is an arbitrary
 chosen orientation of the bounded edges, that is,
 $\partial E = \{\partial^-E, \partial^+E\}$.
For $E\in \Gamma^{[1]}_{\infty}$,
 $\partial^-E$ denotes the unique vertex adjacent to $E$.

The index $\delta_i(\Gamma, \bold{E}, h, \bold{A})$, for each marked edge $E_i$, is
 given by the product
\[
\delta_i(\Gamma, {\bf E}, h, {\bf A}) = 
w_{\Gamma}(E_i)\cdot [\Bbb Zu_{(\partial^-E_i, E_i)}+L(A_i)\cap N: 
(\Bbb Qu_{(\partial^-E_i, E_i)}+L(A_i))\cap N].
\]

\subsubsection{Incidence conditions for 
 holomorphic curves in $X_0$ and $X_1$}\label{subsubsec:incidence}
Next, we define incidence conditions for holomorphic curves.
Assume that $X_0$ is a toric variety defined by a complete fan in 
 $N_{\Bbb Q}$.
We use the same notation as above.
\begin{definition}\label{constraint}
In the case of toric variety $X_0$, we take incidence conditions
 to be the subvarieties of $X_0$
 given as the closures of the orbits of general points 
 $\{q_1, \dots, q_l\}$ in $X_0$,
 by the subtori of the big torus acting on $X_0$
 corresponding to the linear subspace
 $\{L_i\}$.
We write these subvarieties by ${\bf Z} = \{Z_i\}$.
\end{definition}
For $X_1$, recall that there is a gradient Hamiltonian flow
 $\phi_{grH}: X_1\to X_0$.
Since $\phi_{grH}: X_1\to X_0$ is diffeomorphic only away from the singular locus of
 $X_0$, if $Z_i$ intersects the singular locus of $X_0$, the inverse image 
 $\phi_{grH}^{-1}(Z_i)$ may not be of pure dimensional cycle, or may have boundary. 
So we assume the following.\\

\noindent
{\bf Assumption 1.}
Let us write $\dim Z_i = m_i$.
Let $int X_0$ be the complement of the union of toric divisors.
Then the inverse image 
\[
\phi^{-1}_{grH}(Z_i\cap intX_0)
\]
 can be completed to an 
 $m_i$-dimensional cycle in $X_1$ in the following sense.
Let $W = X_0\setminus int X_0$ be the union of toric divisors.
Then there is an $m_i$-dimensional chain $C_i$ in $\phi_{grH}^{-1}(W)$ such that 
 $\phi^{-1}_{grH}(Z_i\cap intX_0)\cup C_i$ is an $m_i$-dimensional
 cycle in $X_1$.

\begin{definition}\label{constraint2} 
In the case of $X_1$, we assume that each $Z_i$ satisfies Assumption 1, 
 and take the cycle $\widetilde Z_i = \phi^{-1}_{grH}(Z_i\cap intX_0)\cup C_i$ as
 an incidence condition.
There are choices of the chain $C_i$,
 and the homology class of $\widetilde Z_i$
 may not be unique. 
We choose one from these choices and write the
 cycle by $\widetilde Z_i$.
\end{definition}
We note that main results in this paper do not depend on the choice of $\widetilde Z_i$. 

When Assumption 1 holds for $Z_i$, 
 then we can assume that there is a continuous family of cycles
 $\mathcal Z_i\subset \mathfrak X$ over $B$ such that:
\begin{itemize} 
\item  $\mathcal Z_i|_{X_1} = \widetilde Z_i$ and $\mathcal Z_i|_{X_0} = Z_i$.
\item The complement in $\mathcal Z_i$ of the union 
\[
\bigcup_{t\in [0, 1]}\phi_{grH, t}^{-1}(Z_i\cap intX_0)
\]
 of the inverse images of
 $Z_i\cap intX_0$ by the maps $\phi_{grH, t}$, $t\in[0, 1]$,
 is contained in 
 the union of the inverse images of the toric divisors 
 $\bigcup_{t\in [0, 1]}\phi_{grH, t}^{-1}(W)$.
\end{itemize}

\subsection{Preliminary arguments about homology classes}
Let $(\mathfrak X, \widetilde{\omega})$
 be a toric Fano degeneration of 
  an $n$-dimensional Fano projective manifold $(X, \omega)$.
Let $X_0$
 be the toric Fano variety.\\
%Let $\{ D_i\}$ be the set of toric divisors on $X_0$.\\

\noindent
{\bf Assumption 2.}
We assume that $X_0$ allows a small resolution.
Also, we assume that the (pointed) base space
 $(B;p_0, p_1)$ of the degeneration $\mathfrak X$
 is an open subset of $\Bbb C^m$ for some $m\in \Bbb N$ which is diffeomorphic
 to the pointed open ball, and the marked points are 
 the origin ($=p_0$) and the point $(1, 0, \dots, 0)=p_1$.
Let $X_t$ be the fiber over $(t, 0, \dots, 0)$.
We assume that $X_t$, $t\in (0, 1]$ is nonsingular.\\

We proved the following results in \cite{NNU}.
\begin{prop}\label{homgr}(Lemma 9.2, \cite{NNU})
\begin{enumerate}[(i)]
\item When $(\mathfrak X, B)$ is a toric degeneration of $X_1$, 
 there is a map $\phi: X_1\to X_0$ which is 
 natural up to homotopy (in particular, it is homotopic to $\phi_{grH}$).
The map $\phi$ is diffeomorphic away from the small neighbourhood of 
 the singular locus of $X_0$. 
\item 
When $X_0$ allows a small resolution, it induces an isomorphism
 $\phi_*: \pi_2(X_1)\to \pi_2(X_0)$.\qed
\end{enumerate}
\end{prop}
Since $X_0$ is a compact toric variety, it is simply connected.
So the natural isomorphism 
\[
\pi_2(X_0)\simeq H_2(X_0, \Bbb Z)
\]
 holds.
By $(ii)$ of Proposition \ref{homgr}, the map $\phi$ induces an epimorphism
 of homology groups
\[
\phi_*: H_2(X_1, \Bbb Z)\to H_2(X_0, \Bbb Z)
\] 
 which has a natural splitting.
Thus, the group $H_2(X_1, \Bbb Z)$ can be written as
\[
H_2(X_1, \Bbb Z)\simeq H_2'\oplus H_2'',
\]
 where $\phi_*|_{H_2'}$ is an isomorphism and
 $\phi_*|_{H_2''}$ is zero.
In particular, the summand $H_2'$ is torsion free.
\begin{rem}
In most of the examples I know, the map $\phi_*$ is an isomorphism.
It might be possible to prove it under some moderate assumption,
 although I do not try to do it in this paper. 
\end{rem}

Since $X_1$ is smooth, Poincar\'e duality holds, 
 so there is a natural isomorphism
\[
(H_{2n-2}(X_1, \Bbb Z))^*\simeq H_2'\oplus fH_2'',
\]
 here $fH_2''$ is the torsion free part of $H_2''$.
Let 
\[
(H_{2n-2}(X_1, \Bbb Z))^*\simeq (H_{2n-2}')^*\oplus (H_{2n-2}'')^*
\]
 be the corresponding splitting of $(H_{2n-2}(X_1, \Bbb Z))^*$.

Let $\pi:\widetilde X_0\to X_0$ be a small resolution.
Let $Z$ be the free abelian group generated by the prime toric divisors of 
 $X_0$.
Since $\pi$ is small, toric prime divisors of $X_0$ and $\widetilde X_0$ 
 are in natural one-to-one correspondence.
So we write by $Z$ the free abelian group generated by the prime toric divisors of 
 $\widetilde X_0$, too.
Thus, there is a commutative diagram
\[\xymatrix{
& H_{2n-2}(\widetilde X_0, \Bbb Z)\ar[d]^{\pi_*}\\
Z\ar[ru]^{p_1}\ar[r]^{p_2\hspace{.3in}} & H_{2n-2}(X_0, \Bbb Z).
}\]
Here $p_1, p_2$ are surjections, and $\pi_*$ is an isomorphism, since $\pi$ is
 small.

Since $\widetilde X_0$ is a toric variety, its second homotopy group is generated by
 holomorphic spheres, and a natural isomorphism
\[
\pi_2(\widetilde X_0)\simeq H_2(\widetilde X_0, \Bbb Z)
\] 
 holds.
In particular, $ H_2(\widetilde X_0, \Bbb Z)$ is free.
Also, since $\widetilde X_0$ is smooth, we have a natural Poincar\'e duality
 isomorphism
\[
(H_{2n-2}(\widetilde X_0, \Bbb Z))^*\simeq H_2(\widetilde X_0, \Bbb Z).
\]
On the other hand, there is a natural surjection
\[
H_2(\widetilde X_0, \Bbb Z)\to H_2(X_0, \Bbb Z).
\]
Combining with $H_2'\simeq H_2(X_0, \Bbb Z)$, 
 we see that there is a natural surjection
\[
H_2(\widetilde X_0, \Bbb Z)\to H_2'.
\]
Taking the Poincar\'e dual of the both sides, we have a surjection
\[
(H_{2n-2}(\widetilde X_0, \Bbb Z))^*\to (H_{2n-2}')^*.
\]
On the other hand, the dual of the map $p_1$ gives an inclusion
\[
p_1^*: (H_{2n-2}(\widetilde X_0, \Bbb Z))^*\to Z^*.
\]
\subsection{Degrees and homology classes} 
First we recall the following definition from \cite{NS}.
\begin{definition}
Let $X_0$ be an $n$-dimensional toric variety.
A holomorphic curve $C\subset X_0$ is \emph{torically transverse} if it is disjoint from 
 all toric strata of codimension greater than one.
A stable map $\varphi: C\to X_0$ is torically transverse if $\varphi^{-1}(int X_0)\subset C$
 is dense and $\varphi(C)\subset X_0$ is a torically transverse curve. 
Here $int X_0$ is the complement of the union of toric divisors.

\end{definition}
Let $\Sigma\subset N_{\Bbb R}$ be the fan defining $X_0$.
For each ray of $\Sigma$, we have the generator $v\in N$ and 
 its associated toric prime divisor $D_v$.
\begin{definition}\label{def:degforvarphi}
For a torically transverse curve $\varphi: C\to X_0$,
 the \emph{degree} is given by the map
 \[\Delta(\varphi): N\setminus\{0\}\to\Bbb N\] defined as follows.
For a primitive $v\in N$ and $\lambda\in \Bbb N$, 
 $\lambda v$ is mapped to 0 if $\Bbb R_{\geq 0}v$ is not a ray of $\Sigma$,
 and to the number of points of multiplicity $\lambda$ in $\varphi^*D_v$ otherwise.
\end{definition}
Recall the degree \[\Delta: N\setminus\{0\}\to \Bbb N\]
 of a tropical curve in $\Bbb R^n$ was given by the data of the direction 
 vectors and multiplicity of the unbounded edges.
We have to define the notion of degree for curves in a general fiber $X_t$.
\begin{definition}\label{deg_t}
We define the degree of a 
 curve in general smooth fiber $X_t$, $t\neq 0$ to be its 
 integral homology class.
\end{definition}

Let $\varphi: C\to X_0$ be a torically transverse 
  stable map of degree $\Delta$.
We define a map $\Delta_D$ from the set of primitive vectors in
 $N\setminus\{0\}$ to $\Bbb N$ as follows.
Namely, to a primitive vector $v\in N$,
 we associate
\[
\Delta_D(v) = \sum_{a>0}a\Delta(av) \in \Bbb N.
\]
This can be regarded as an element of the dual space $Z^*$ of 
 the space $Z$ introduced above.

\begin{definition}
A \emph{coarse-degree} 
 is a map $\widetilde{\Delta}$ from the set of
 nonzero primitive vectors in $N$ to
 the set of non-negative integers satisfying
 the condition
 \[\sum_{v:primitive}\widetilde{\Delta}(v)v = 0.\]
This can be extended to an element of $Z^*$, and we write it by the same letter
 $\widetilde\Delta$.
\end{definition}
A degree $\Delta$
 of a tropical curve or a 
 torically transverse stable map in a toric variety
 naturally gives a coarse-degree:
\begin{lem}
For a degree $\Delta$, the map $\Delta_{D}$ is a coarse-degree.
\end{lem}
\proof
This follows from the balancing condition for the tropical curve.\qed
\begin{definition}
A coarse-degree $\widetilde{\Delta}$ is called \emph{rational}
 if it is induced from a degree of a rational tropical curve (or of a torically transverse
 rational stable map in some toric variety).
Let $\mathfrak D\subset Z^*$ be the set of rational coarse-degrees of $X_0$.
\end{definition}
Recall that there is an inclusion
\[
p_1^*: (H_{2n-2}(\widetilde X_0, \Bbb Z))^*\to Z^*.
\]
\begin{lem}\label{deguni}
The set $\mathfrak D\subset Z^*$
 is a subset of (the image of) the space $(H_{2n-2}(\widetilde X_0, \Bbb Z))^*$.
\end{lem}
\proof
Let 
\[
p_1: Z\to H_{2n-2}(X_0;\Bbb Z)
\]
 be the quotient map defined before.
Given an element $\widetilde{\Delta}$
 of $\mathfrak D$, considered as a linear function on
 $Z$, 
 define 
$\widetilde{\Delta}'$ by 
\[\widetilde{\Delta}'(p_1(D)) = 
 \widetilde{\Delta}(D).\]
This defines a well-defined element of 
 $H_{2n-2}(X_0;\Bbb Z)^{*}$, because if $p_1(D)$ is zero,
 $\widetilde{\Delta}(D)$ have to be zero because 
 $\widetilde{\Delta}(D)$ is the sum of the 
 transversal intersection numbers of 
 a rational curve and a linear sum $D$ of toric divisors, 
 so it must be zero when $D$ is homologous to zero.\qed
\begin{lem}
$\mathfrak D$ is a submonoid of $(H_{2n-2}(\widetilde X_0;\Bbb Z))^*$.
\end{lem}
\proof
It suffices that the set of degrees of rational tropical curves forms a monoid.
This follows because given two rational tropical curves, we can parallel transport 
 one of them so that two tropical curves intersect.
By taking a suitable union of graphs as a domain,
 we have a new rational tropical curve whose degree is the sum of the given two.\qed\\

On the other hand,
 by the Fredholm regularity of rational stable maps in a toric variety
 \cite{CO} (there, the Fredholm regularity for
 holomorphic disks with Lagrangian boundary condition in toric manifold is proved.
 The rational curve case follows from more simple
 cohomological argument, and since
 a torically transverse curve is contained
 in the smooth part of $X_0$, we can use the regularity result),
 we can deform torically transverse rational stable maps in $X_0$ into
 rational stable maps in $X_1$
 (precisely speaking, to use the regularity result, we should take $X_t$ with small $t$
 instead of $X_1$.
 Reparameterizing the base space
 of the degeneration if necessary, we assume $X_1$ is sufficiently close to $X_0$).

A torically transverse rational stable map $\varphi: C\to X_0$ lifts to 
 a torically transverse rational stable map $\widetilde{\varphi}: C\to \widetilde X_0$.
Since $\widetilde X_0$ is smooth, 
 the image of $\widetilde{\varphi}$ gives an element 
 of $(H_{2n-2}(\widetilde X_0, \Bbb Z))^*$.
On the other hand, by the above argument, $\varphi$ lifts to a
 torically transverse stable map $\varphi_1: C\to X_1$.
Since $X_1$ is smooth, too, $\varphi_1$ gives an element of 
 $(H_{2n-2}(X_1, \Bbb Z))^*$,
 in fact, an element of $(H_{2n-2}')^*$ using the notation of the previous subsection.

Recall that there is a natural surjection
\[
P: (H_{2n-2}(\widetilde X_0, \Bbb Z))^*\to (H_{2n-2}')^*
\]
 and a submonoid
\[
\mathfrak D\subset (H_{2n-2}(\widetilde X_0, \Bbb Z))^*
\]
Using the construction so far, we see the following.
\begin{prop}\label{degree}
Let $[\varphi']\in H_2'$ be an element represented by a rational stable map
 which is a lift of some torically transverse
 rational stable map $\varphi$ in $X_0$.
Then the intersection of $P^{-1}([\varphi'])$ and $\mathfrak D$ in 
 $(H_{2n-2}(\widetilde X_0, \Bbb Z))^*$, or its image in $Z^*$ is
 given as follows:
\begin{enumerate}
\item Consider all the torically transverse
 rational stable maps 
\[
\varphi_i: C_i\to X_0, \;\;i = 1, \dots, m 
\]
 whose lifts to $X_1$ give the class $[\varphi']$.
\item Lift $\varphi_i$ to $\widetilde X_0$.
It gives an element $a_i$ of $Z^*$.
\item Then we have
\[
P^{-1}([\varphi'])\cap \mathfrak D = \{a_1, \dots, a_m\}.
\]
\end{enumerate}
\qed
\end{prop}
\begin{definition}\label{def:degree}
For $[\varphi']\in H_2'$ as in the proposition above, we call the subset
\[
P^{-1}([\varphi'])\cap\mathfrak D
\]
 as the set of \emph{coarse-degrees corresponding to \text{$[\varphi']$}}.
We may also consider it as a subset of $Z^*$ by the natural inclusion $p_1^*$. 
\end{definition}

\subsection{Stable maps in $X_0$ and their deformations to $X_1$}

Let us assume we have fixed a degree $\Delta$ of a torically transverse
 rational stable map 
  in $X_0$.
Let $e=|\Delta|$ be the number of the unbounded edges of the associated
 tropical curve as before.
Take generic incidence conditions $(Z_1, \dots, Z_l)$ in $X_0$
 corresponding to a sequence of linear subspaces $(L_1, \dots, L_l)$ of $N_{\Bbb R}$.
Write $\text{codim}\, L_i = \text{codim} \, Z_i = d_i+1$
 and assume $\sum_{i=1}^l d_i = e+n-1$.
In the following lemma, we do not a priori assume that $\varphi$ is torically transverse.
\begin{lem}\label{transverse}%KEY LEMMA
Let $\varphi: C\to X_0$ be a rational stable map of degree $\Delta$
 satisfying the 
 incidence conditions $(Z_1, \dots, Z_l)$.
Assume that each $Z_i$ is transversal to each toric stratum.
Then the domain curve $C$ is the nonsingular rational 
 curve, and the map $\varphi$ is torically transverse.
\end{lem}
\proof
Let $\varphi: C\to X_0$ be any (not necessarily torically transverse) 
 rational stable map satisfying the incidence conditions $(Z_1 , . . . , Z_l)$,
 and assume that the homology class of the image $[\varphi(C)]$
 is the same as the homology class of the torically transverse rational stable map
 of degree $\Delta$.
Then by the assumption that $Z_i$ is transversal to the toric strata and
 that $X_0$ is Fano, 
 one sees that the domain curve $C$ is the nonsingular rational curve by
 dimension counting.

Now suppose that the image $\varphi(C)$ is contained in some lower dimensional
 toric stratum $D$.
Then each $Z_i$ must intersect $D$, and by the transversality assumption, 
 it gives a codimension $d_i+1$ condition in $D$.
On the other hand, the map $\varphi$, considered as a map to $D$, 
 has moduli space whose virtual dimension is lower than 
 $\varphi$ considered as a map to $X_0$.
Thus, such a map cannot satisfy the generic incidence conditions.
So the image $\varphi(C)$ intersects the toric divisors in discrete points.

In \cite{NNU}, it was proved that under Assumption 2, 
 a stable map with this property whose image intersects lower dimensional 
 toric strata belongs to a lower dimensional subvariety in the moduli,
 so they cannot satisfy generic incidence conditions.
Thus, $\varphi$ must be torically transverse.\qed
\begin{cor}
When each variety $Z_i$ is a point, then any rational stable map $\varphi: C\to X_0$
 satisfies the properties that $\varphi$ is torically transverse and the 
 domain curve $C$ is the nonsingular rational curve.
\end{cor}
\proof
When each $Z_i$ is a point, the transversality assumption is equivalent to
 the statement  that $Z_i$ is not contained in the union of toric divisors.
However, this is automatically satisfied when $\{Z_i\}$ is generic.\qed\\

As we mentioned before, 
 by the Fredholm regularity of torically transverse 
 stable maps in a toric variety
 (\cite{CO}), we can deform holomorphic curves in $X_0$ into
 $X_1$, and this correspondence is (non canonically)
 one to one.
In particular, we have the following.
Note that we assume $X_1$ is sufficiently close to $X_0$, 
 as we remarked in the previous subsection.
\begin{prop}\label{lift}
Let $\varphi: C\to X_0$ be a stable map of given degree $\Delta$ as above,
 satisfying the 
 incidence conditions $(Z_1, \dots, Z_l)$,
 and assume that each $Z_i$ is transversal to each toric stratum,
 as Lemma \ref{transverse}.
Assume also that
 each $Z_i$ satisfies Assumption 1 and choose $\widetilde Z_i\subset X_1$
 for each $i$.
Then
 $\varphi$ can be uniquely deformed into a stable map $\varphi_1: C\to X_1$ 
 satisfying the incidence conditions $(\widetilde Z_1, \dots, \widetilde Z_l)$.
Its degree $\beta$ satisfies the property that the coarse-degree $\Delta_D\in Z^*$
 associated to $\Delta$ is contained in the set of 
 coarse-degrees corresponding to $\beta$ (Definition \ref{def:degree}).
\qed
\end{prop} 
For each $i$, we take a family of cycles $\mathcal Z_i$ over the base $B$
 as in Subsection \ref{subsubsec:incidence}.
Then the above proposition can be formulated for this family, 
 giving lifts $\varphi_t: C\to X_t$ for every $t\in [0, 1]$.
This does not depend on the choice of $\mathcal Z_i$.
We can prove the converse of Proposition \ref{lift}.
\begin{prop}\label{prop:converse}
Consider the same situation as Proposition \ref{lift}.
Then there is a positive number $T$ with the following property:
For any $t$ with $0<t<T$, assume that there is a rational stable map 
 $\varphi_t: C\to X_t$ of degree $\beta$ satisfying the incidence conditions
 given by the restrictions $\mathcal Z_i|_{X_t}$. 
Then $\varphi_t$ is contained in the ones constructed in Proposition \ref{lift}
 (with $X_1$ replaced by $X_t$).
\end{prop}
\proof
Suppose the statement is false.
Then there is a sequence of rational stable maps 
 $\psi_{t_i}: C_i\to X_{t_i}$ of degree $\beta$, with $t_i\to 0$ as
 $i\to\infty$
 satisfying the incidence conditions,
 but not contained in the ones constructed in Proposition \ref{lift}.
By Gromov's compactness theorem, we may assume that $\psi_{t_i}$
 converges to a limit rational stable map $\psi: C'\to X_0$ 
 satisfying the incidence conditions $(Z_1, \dots, Z_l)$,
 whose homology class is $\phi_*(\beta)$.
 
The virtual dimension of a rational stable map in $X_0$,
 whose homology class is $\phi_*(\beta)$, is $e+n-1$,
 which is the same as the rational stable maps of degree $\Delta$. 
Then as in the proof of Lemma \ref{transverse}, $\psi$ 
 must be torically transverse and $C' = \Bbb P^1$.
Then since the lift of $\psi$ to $X_t$ is unique by Proposition \ref{lift}, 
 $\psi_{t_i}$ must be contained in the family constructed in Proposition \ref{lift}.\qed

\subsection{Counting invariants via toric Fano degeneration}
We define two counting invariants.
One is directly related to the genuine Gromov-Witten invariants, and the other is 
 the transversal Gromov-Witten type invariants considered in \cite{NS}.

\subsubsection{Gromov-Witten invariants}
As before, reparameterizing the base
 $B$, we assume $X_1$ satisfies the condition of
 Proposition \ref{prop:converse}.

We briefly recall the definition of Gromov-Witten invariants.
See for example \cite{CK}, for more information.
Let $X$ be a projective algebraic manifold and $\beta\in H_2(X, \Bbb Z)$.
Let $\overline{\mathcal M}_{g, n+1}(X, \beta)$ be the moduli stack
 of stable maps of genus $g$ with target $X$, the homology class of the 
 image is $\beta$, 
 and let $[\overline{\mathcal M}_{g, n}(X, \beta)]^{virt}$ be its virtual fundamental class.
$[\overline{\mathcal M}_{g, n}(X, \beta)]^{virt}$ has dimension
 \[(1-g)(\dim X-3)-\int_{\beta}\omega_X + n.\]
Let 
\[\pi_1: \overline{\mathcal M}_{g, n}(\beta, X)\to X^n\]
 be the natural map given by evaluations.
Let $\alpha_1, \dots, \alpha_n$ be elements of $H^*(X, \Bbb Q)$.
Then the Gromov-Witten invariant 
 $\langle I_{g, n, \beta}\rangle(\alpha_1, \dots, \alpha_n)$ is defined by
\[
 \langle I_{g, n, \beta}\rangle(\alpha_1, \dots, \alpha_n)
   = \int_{[\overline{\mathcal M}_{g, n}(X, \beta)]^{virt}}
    \pi_1^*(\alpha_1\otimes\cdots\otimes \alpha_n).
  \]
We are concerned with the case when $g = 0$.
As above, a torically transverse rational stable map
 $\varphi$ in $X_0$ lifts to $\varphi_1$ in $X_1$.
The following is almost clear.
\begin{lem}\label{virdim}
The virtual dimension of $\varphi$ is the same as that of $\varphi_1$.
\qed
\end{lem}

To state the main theorem, we recall what we count
 and at the same time we introduce some notations.
Let us fix a class $\beta\in H_2'$ represented by a rational stable map.
Let $\Delta_{\beta} = 
 \{a_1, \dots, a_m\}$ be the set of coarse-degrees corresponding to $\beta$
 in the sense of Definition \ref{def:degree}.
 
On $X_0$, we count (a priori not necessarily torically transverse)
 stable maps
 $\varphi: C\to X_0$
 such that:
\begin{enumerate}
\item The domain $C$ is a connected, $l$-pointed rational prestable curve.
\item $\varphi$ satisfies general incidence conditions $\bold Z = (Z_1, \dots, Z_l)$
 of codimension ${\bold d} = (d_1, \dots, d_l)$, 
 satisfying Assumption 1 and the transversality assumption of Lemma \ref{transverse}.
\item Then by Lemma \ref{transverse}, $\varphi$ has to be torically transverse.
In particular, the degree $\Delta$
 of $\varphi$ in the sense of Definition \ref{def:degforvarphi}
 is defined.
The degree $\Delta$ should satisfy the property that
 the associated coarse-degree $\Delta_D$ is contained in the set
 $\Delta_{\beta}$.
\end{enumerate}
By Lemma \ref{transverse},
 $C$ has to be the nonsingular rational curve.

For the case of $X_1$, we count stable maps $\psi: C\to X_1$ such that:
\begin{enumerate}
\item $C$ is a connected, $l$-pointed rational prestable curve.
\item  $\psi$ 
  satisfies the incidence conditions determined by Definition \ref{constraint2} 
  (which we write by $\bold{\widetilde Z}$).
\item $\psi$
   has degree $\beta\in H_2'\subset H_2(X_1, \Bbb Z)$.
\end{enumerate}
For the case of tropical curves, we count tropical curves $h: \Gamma\to N_{\Bbb R}$
 such that
\begin{enumerate}
\item The graph $\Gamma$ is connected, rational, and
 $l$-marked: ${\bf E} = (E_1, \dots, E_l)$.
\item $h: \Gamma\to\Bbb R^n$
  satisfies generic incidence conditions ${\bf A} = (A_1, \dots, A_l)$
  of codimension ${\bf d}$, such that $Z_i$ is the closure of an orbit
  of the subtorus corresponding to $L_i$, the linear subspace of $N_{\Bbb R}$
  parallel to $A_i$.
\item $(\Gamma, h)$ has degree $\Delta$ and this $\Delta$ satisfies the 
 same property as in the case of $X_0$.
\end{enumerate}
Tropical curves have to be counted with multiplicity 
\[
w(\Gamma, \bold E)\cdot
 \mathfrak D(\Gamma, \bold{E}, h, \bold{A})\cdot\prod_{i=1}^l\delta_i
 (\Gamma, \bold{E}, h, \bold{A})
\]
 (see Subsection \ref{subsubsec:trop}).
 
We write these three counting numbers by
 $N_{X_0}(\Delta_{\beta}, \bold Z)$, $N_{X_1}(\beta, \bold{\widetilde Z})$
 and $N_{trop}(\Delta_{\beta}, \bold L)$
respectively.
If the virtual dimension of the moduli spaces of the
 curves is not equal to the codimension of 
 the incidence conditions, we define these numbers to be 0.
Our main theorem is the following.
\begin{thm}\label{tropical}
The following equalities hold.
\[N_{X_1}(\beta, \bold Z) = N_{X_0}(\Delta_{\beta}, \bold{\widetilde Z})
  = N_{trop}(\Delta_{\beta}, \bold A).
  \]
These numbers do not depend on the choice of $\bold Z$,
 $\bold{\widetilde Z}$ or $\bold A$,
 if they are general.
\end{thm}
\proof
In \cite{NS}, Theorem 8.3 and Corollary 8.4,
 the equality between $N_{X_0}(\Delta_{\beta}, \bold Z)$
 and $N_{trop}(\Delta_{\beta}, \bold A)$,
 as well as the independence to the choice of incidence conditions are shown. 
So it suffices to show the equality between 
 $N_{X_1}(\beta, \bold{\widetilde Z})$ and 
 $N_{X_0}(\Delta_{\beta}, \bold Z)$.
 
By Proposition \ref{lift},
 we see that any rational stable map of degree $\Delta$
 such that $\Delta_D$ is contained in $\Delta_{\beta}$ and 
 satisfying the incidence conditions ${\bf Z}$ in $X_0$ uniquely lifts
 to a rational stable map 
 in $X_1$ satisfying the incidence conditions ${\bf{\widetilde Z}}$.
By Lemma \ref{virdim}, the virtual dimensions of 
 the rational stable map in $X_0$ and its lift in $X_1$ are the same.
This shows the inequality
 \[N_{X_1}(\beta, \bold{\widetilde Z})\geq 
  N_{X_0}(\Delta_{\beta}, \bold Z)\] holds.
So the problem is to show that this lift is surjective.
But it is the content of Proposition \ref{prop:converse}.\qed
\begin{cor}\label{GW}
The number $N_{X_1}(\beta, \bold{\widetilde Z})$ is the Gromov-Witten invariant
 $I_{0, l, \beta}(\alpha_1, \dots, \alpha_l)$, where $\alpha_i$ is the Poincar\'e  dual
 of the homology class of ${{\widetilde Z_i}}$.
\end{cor}
\proof
By the theorem, the stable maps contributing to 
 $N_{X_0}(\Delta_{\beta}, \bold Z)$ and $N_{X_1}(\beta, \bold{\widetilde Z})$
  are in
 one-to-one correspondence.
Recall that the stable maps contributing to $N_{X_0}(\Delta_{\beta}, \bold Z)$
 are Fredholm regular.
So the stable maps contributing to $N_{X_1}(\beta, \bold{\widetilde Z})$,
 which are small perturbations of those contributing to 
 $N_{X_0}(\Delta_{\beta}, \bold Z)$,
 are Fredholm regular, too
 (recall that we assume $X_1$ is sufficiently close to $X_0$).
So the moduli space of rational stable maps of degree $\beta$ 
 intersecting with general incidence conditions $\bold{\widetilde Z}$ in $X_1$ is 
 $N_{X_1}(\beta, \bold{\widetilde Z})$ points, and its virtual fundamental cycle is
 equal to itself.\qed\\
\begin{thm}
Using the same notation as in Corollary \ref{GW}, 
The Gromov-Witten invariant $I_{0, l, \beta}(\alpha_1, \dots, \alpha_l)$ for
 $X = X_1$ is equal to  $N_{trop}(\Delta_{\beta}, \bold A)$.
\end{thm}
\proof
This follows from Theorem \ref{tropical} and Corollary \ref{GW}.\qed
%%%%%%%%%%%%%%%%%%%%%%%%%%%%%%%%%%%%%%%%%%
\subsubsection{Transversal Gromov-Witten type invariants}
The transversality assumption of Lemma \ref{transverse} for an incidence condition
 $Z_i$ in $X_0$ does not hold in general.
Nevertheless, even without this assumption we can define counting invariants
 of curves by restricting attention only to torically transversal curves,
 as we did in \cite{NS}.
When the transversality assumption holds, all the curves satisfying 
 the incidence conditions have to be torically transverse, so 
 under this assumption, these two invariants are equal.

Thus, here we only assume Assumption 1 for the incidence conditions $Z_i$.
In this case, what we count are the following.
As before, we fix a class $\beta\in H_2'$ represented by a rational stable map.

On $X_0$, we count 
 stable maps
 $\varphi: C\to X_0$
 such that:
\begin{enumerate}
\item The domain $C$ is a connected, $l$-pointed rational prestable curve.
\item $\varphi$ is torically transverse
 and satisfies general incidence conditions $\bold Z = (Z_1, \dots, Z_l)$
 of codimension ${\bold d} = (d_1, \dots, d_l)$, 
 satisfying Assumption 1.
\item $\varphi$ has a fixed degree $\Delta$
 such that the associated coarse-degree $\Delta_D$ is contained in the set 
 $\Delta_{\beta}$.
\end{enumerate}
By Lemma \ref{transverse},
 $C$ has to be the nonsingular rational curve.
 
For the case of $X_1$, we count stable maps $\psi: C\to X_1$ such that:
\begin{enumerate}
\item $C$ is a connected, $l$-pointed rational prestable curve.
\item  $\psi$ 
  satisfies the incidence conditions determined by Definition \ref{constraint2} 
  (which we write by $\bold{\widetilde Z}$)
 and $\psi$ is a lift of a torically transverse stable map $\varphi$ in $X_0$
   as in Proposition \ref{lift}.
   Note that in Proposition \ref{lift}, the transversality assumption for $Z_i$ is assumed,
   but for torically transverse $\varphi$, the conclusion holds without this assumption.
\item $\psi$
   has degree $\beta\in H_2'\subset H_2(X_1, \Bbb Z)$.
\end{enumerate}
For the case of tropical curves, we count the same
 objects as before, namely, tropical curves $h: \Gamma\to N_{\Bbb R}$
 such that
\begin{enumerate}
\item The graph $\Gamma$ is connected, rational, and
 $l$-marked: ${\bf E} = (E_1, \dots, E_l)$.
\item $h: \Gamma\to\Bbb R^n$
  satisfies generic incidence conditions ${\bf A} = (A_1, \dots, A_l)$
  of codimension ${\bf d}$, such that $Z_i$ is the closure of an orbit
  of the subtorus corresponding to $L_i$, the linear subspace of $N_{\Bbb R}$
  parallel to $A_i$.
\item $(\Gamma, h)$ has degree $\Delta$ and this $\Delta$ satisfies the 
 same property as in the case of $X_0$.
\end{enumerate}
As we mentioned before, 
 tropical curves have to be counted with multiplicity.
 
We write these three counting numbers by
 $N^{trans}_{X_0}(\Delta, \bold Z)$, $N^{trans}_{X_1}(\beta, \bold{\widetilde Z})$
 and $N_{trop}(\Delta, \bold L)$
respectively.
Then the conclusion as Theorem \ref{tropical} is proved by the same proof.
\begin{thm}\label{tropical2}
The following equalities hold:
\[N^{trans}_{X_1}(\beta, \bold Z) 
 = N^{trans}_{X_0}(\Delta_{\beta}, \bold{\widetilde Z})
  = N_{trop}(\Delta_{\beta}, \bold A).
  \]
These numbers do not depend on the choice of $\bold Z$,
 $\bold{\widetilde Z}$ or $\bold A$,
 if they are general.\qed
\end{thm}
\begin{rem}\label{rem:GW}
\begin{enumerate}
\item 
Note that $N^{trans}_{X_1}(\beta, \bold Z)$ may not be a homological invariant, 
 that is, if we change the cycle $\widetilde Z_i$ within the same homology class,
 but not necessarily related to $Z_i$
 as in Assumption 1, the number of rational stable maps incident
 to them may change.
However, by the Fredholm regularity of the maps contributing to 
 $N^{trans}_{X_1}(\beta, \bold Z)$, 
 it gives a lower bound for the Gromov-Witten invariants:
\[
N^{trans}_{X_1}(\beta, \bold Z)\leq I_{0, l, \beta}(\alpha_1, \dots, \alpha_l),
\]
 where $\alpha_i$ is the Poincar\'e dual of the homology class of $\widetilde Z_i$.
\item An interesting case to apply this construction is the (full) flag manifolds of type A.
According to \cite{KM}, Schubert subvarieties degenerates to unions of toric
 strata under toric degeneration.
Since the general linear group acts on the set of Schubert subvarieties
 by multiplication, the degenerated one (a union of toric strata) may also be deformed, 
 and it may be possible to represent the homology class of it by a union of 
 the classes $[Z_i]$, where $Z_i$ is the closure of an orbit 
 the action of a torus as above.
Then the problem is whether the number $N^{trans}_{X_1}(\beta, \bold Z)$
 is genuinely smaller than $I_{0, l, \beta}(\alpha_1, \dots, \alpha_l)$ or not, 
 where $\alpha_i$ is the Poincar\'e dual of the corresponding Schubert subvariety.

On the other hand, some of the Schubert subvarieties are homologous to
 the cycles of the form $\widetilde Z_i$, where $Z_i$ is a cycle in $X_0$ 
 which satisfies Assumption 1 and the transversality assumption of Lemma
 \ref{transverse}. 
For these cases, we can calculate the genuine Gromov-Witten invariants
 (see the next section).
\end{enumerate}
\end{rem}

\section{Examples}
\subsection{Flag manifold $F_3$}
In this section, we give an example of calculation
 of Gromov-Witten invariants by tropical method in the case of flag manifold
 $F_{3}$, which parameterizes the full flags in $\Bbb C^3$.
It is embedded in $\Bbb P^2\times \Bbb P^2$
 by Pl\"ucker embedding, 
 and its toric degeneration is given explicitly by
\[
  \mathfrak{X} =
  \Bigl\{ 
    \bigl([Z_1 : Z_2 : Z_3], [Z_{12} : Z_{13} : Z_{23}] ,t \bigr)
   \, \Bigm| Z_1 Z_{23} - Z_2 Z_{13} + t Z_3 Z_{12} =0 \,
   \Bigr\}.
\]
with the central fiber
\[
  X_0 = \Bigl\{ 
    \bigl([Z_1 : Z_2 : Z_3], [Z_{12} : Z_{13} : Z_{23}] \bigr)
    \in \mathbb{P}^2 \times \mathbb{P}^2 
   \, \Bigm|  Z_1 Z_{23} - Z_2 Z_{13}  = 0 \,
   \Bigr\}.
\]
The Gelfand-Cetlin polytope (see \cite{NNU} for details) of $F_3$ is defined by the inequalities
\[
  \begin{matrix}
    \lambda_1 && > && 
    \hspace{-10pt} \lambda_2 && > && 
    \hspace{-10pt} \lambda_3 \\
    & {\hspace{-10pt} \uge} && \hspace{-10pt} \dge 
    && \hspace{-10pt} \uge  && 
      {\hspace{-10pt} \dge}&\\
    && \hspace{-10pt} \hbox to1.05em{$\lambda_1^{(2)}$} &&&&
    \hspace{-10pt} \hbox to1.05em{$\lambda_2^{(2)}$}  && \\
    &&& \hspace{-10pt} \uge && \hspace{-10pt} \dge &&& \\
  &&&& \hspace{-10pt} \hbox to1.05em{$\lambda_1^{(1)}$} &&&& 
  \end{matrix}.
\]
and graphically given by the following figure.
\begin{figure}[h]
 \begin{center}
  \includegraphics*{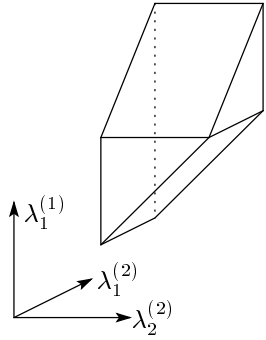}
 \end{center}
 \caption{The Gelfand-Cetlin polytope for $n=3$}
 \label{GCpoly}
\end{figure}
This is the same as the polytope associated to the toric Fano variety
 $X_0$.

We use the coordinate
 $(x, y, z) = (\lambda_2^{(2)}, \lambda_1^{(2)}, \lambda_1^{(1)})$
 for $M_{\Bbb R}\simeq\Bbb R^3$.
Let $\pi_1$, $\pi_2$ be the natural projections of $F_3$ to
 the factors of $\Bbb P^2\times \Bbb P^2$.
Generic fiber of $\pi_i$, $i=1, 2$, is isomorphic to $\Bbb P^1$
 and the corresponding homology classes generate 
 the second homology group $H_2(F_3, \Bbb Z)$ which is isomorphic to $\Bbb Z^2$.
We write these generators by $(1, 0)$ and $(0, 1)$.
Then clearly 
\[I_{0, 1, (1, 0)}(pt) = I_{0, 1, (0, 1)}(pt) = 1.\]
The two point function $I_{0, 1, (1, 1)}(pt, pt)$
 is calculated by counting tropical curves with degree
\[(-1, 0, 0), (0, 1, 0), (0, -1, 1), (1, 0, -1).\]
These are the same as the tropical curves corresponding to the lines in $\Bbb P^3$,
 and so there is one such tropical curve through generic two points.
Thus,  
\[I_{0, 1, (1, 1)}(pt, pt)=1.\]

We calculate the three point functions $I_{0, 3, \beta}(pt, pt, pt)$
 for $\beta = (1, 2), (2, 1)$.
It is easy to see that the tropical curves corresponding to $(1, 2)$
 has degree 
 \[(-1, 0, 0), (0, -1, 1), (0, -1, 0), (0, 1, 0), (0, 1, 0), (1, 0, -1).\]
Projecting such a tropical curve $\bar{\Gamma}$
 to the $xz$-plane, we still have a tropical curve $\Gamma$.
It is a plane tropical curve with degree
 \[(-1, 0), (0, 1), (1, -1).\]
That is, it is a tropical line with unbounded edges of these directions.
Three unbounded edges (directions $(0, -1, 0), (0, 1, 0), (0, 1, 0)$)
 are projected to points on the $xz$-plane.

Projecting generic three points $\bar a, \bar b, \bar c$
 in $\Bbb R^3$ to this plane, 
 we have generic three points $a, b, c$ on $\Bbb R^2$.
If the original tropical curve intersects generic three points, 
 the projected curve also intersects three points.
However, there is no tropical line on $\Bbb R^2$ which is incident
 to generic three points, 
 so it follows that there is no rational tropical curve in $\Bbb R^3$of 
 the degree given above which is incident to $\bar a, \bar b, \bar c$.
This means
\[I_{0, 3, (1, 2)}(pt, pt, pt)=0.\]
The identity 
\[I_{0, 3, (2, 1)}(pt, pt, pt)=0\]
 can be shown in the same way.
By the same proof, it follows that
 \[I_{0, k, (1, k-1)}(pt, \dots, pt) =
 I_{0, k, (k-1, 1)}(pt, \dots, pt) = 0\] for all $k>2$.
 
In fact, we can calculate all genus zero Gromov-Witten invariants
 by tropical method.
As noted above, the second homology group of the flag manifold $F_{3}$ is rank two, 
 and it is generated by rational curves whose classes we wrote by 
 $(1, 0)$ and $(0, 1)$.
We rewrite these by $l_1, l_2$.

In $X_0$, these classes are tropically represented by lines with direction vectors
 $(1, 0, 0)$ and $(0, 1, 0)$.
These curves lift to $X_1$, so, considered as incidence conditions,
 they satisfy Assumption 1 
 and the transversality condition of Lemma \ref{transverse}. 
Real four dimensional classes are evaluated using the divisor axiom.
Generators of these classes $s_1, s_2$ are uniquely determined by the relations
 \[s_i\cap l_j = \delta_{ij}.\]
 
For a class $\beta = (s, t)\in H_2(F_3, \Bbb Z)$, 
 the set of coarse-degrees $\Delta_{\beta}$
 corresponding to $\beta$ (Definition \ref{def:degree}) is given by the following set of 
 degrees for rational tropical curves:
\begin{lem}
A coarse-degree $\widetilde{\Delta}$ belongs to $\Delta_{\beta} = \Delta_{(s, t)}$
 if and only if the following conditions are satisfied:
 
Representing a coarse-degree as an integer valued function on the set of 
 primitive vectors in $N = \Bbb Z^3$, 
\[\begin{array}{ll}
\widetilde{\Delta}(-1, 0, 0) = s, & \widetilde{\Delta} (0, 1, 0) = t,\\
\widetilde{\Delta}(1, 0, 0) = s_1, & \widetilde{\Delta} (1, 0, -1) = s_2,\\
\widetilde{\Delta}(0, -1, 0) = t_1, & \widetilde{\Delta} (0, -1, 1) = t_2,
\end{array}
\]
 where $s_1+s_2 = s$, $t_1+t_2 = t$ and $t_2 = s_2$.
The others are zero.\qed
\end{lem}
Then we have the following.
\begin{prop}
Let 
\[pt\in H_0(F_3, \Bbb Z), l_1, l_2\in H_2(F_3, \Bbb Z), s_1, s_2\in H_4(F_3, \Bbb Z)\]
 be the generators as above
 and 
\[pt^{\vee}, l_1^{\vee}, l_2^{\vee}, s_1^{\vee}, s_2^{\vee}\] be their dual
 cohomology classes.
 
Let $\beta = (s, t)$ be the degree of the curve as above.
Then, when 
\[2(s+t) = 2a+b+c,\] 
 we have,
 \[I_{0, k, \beta}((pt^{\vee})^a, (l_1^{\vee})^b,
  (l_2^{\vee})^c, (s_1^{\vee})^d, (s_2^{\vee})^e)
 = s^dt^e
    I_{0, k', \beta}((pt^{\vee})^a, (l_1^{\vee})^b,
  (l_2^{\vee})^c),\] 
  where 
\[
k = a+b+c+d+e,\;\; k' = a+b+c.
\]

The number
 \[I_{0, k', \beta}((pt^{\vee})^a, (l_1^{\vee})^b,
  (l_2^{\vee})^c)\]
  is the same as the number of $k'$-marked rational tropical curves
  of degree $\Delta$  
  incident to generic $a$ points, generic $b$ lines with direction
  $(1, 0, 0)$ and generic $c$ lines with direction $(0, 1, 0)$.
The degree $\Delta$ should satisfy the condition that the associated coarse-degree
  $\Delta_D$ belongs to $\Delta_{(s, t)}$.
  When $2(s+t)\neq 2a+b+c$, then the invariant is zero.\qed
\end{prop}

\subsection{
Moduli space of rank 2 bundles on a curve of genus two
}
In this subsection, we give another example of calculation of Gromov-Witten invariants
 of a Fano manifold of particular interest
  via tropical method.
Newstead \cite{N}
and Narasimhan and Ramanan \cite{NR}
showed that the moduli space of stable rank two vector bundles
with a fixed determinant of odd degree
on a genus two curve defined as the double cover of $\bP^1$
branched over 
\[\{ \omega_0, \dots, \omega_5 \} \subset \bC,\]
is a Fano complete intersection 
\[X = Q_1 \cap Q_2\] of two quadrics
$$
 Q_1 : \sum_{i=0}^{5} x_i^2 = 0,\hspace{.2in}
 Q_2 : \sum_{i=0}^{5} \omega_i x_i^2 = 0
$$
in $\bP^5$.
The Betti numbers of $X$ are easy to calculate:
\[
b_0=b_2=b_4=b_6=1,\;\; b_3=4,
\]
 and the others are 0.
$X$ has several different structures of integrable systems, known as Goldman systems
 \cite{G}, and $X$ can be degenerated (without considering integrable system)
 to Fano toric varieties in several ways.
Nevertheless, 
 there seems to be no known toric degeneration of $X$ as a Fano integrable system
 in the sense of \cite{NNU}.
In fact, we see that from the view point of toric degenerations of integrable systems
 and enumerations of holomorphic curves associated to them,
 Goldman systems are not good ones (see below).

By deforming the defining equations $Q_1$ and $Q_2$ (and
 doing  changes of variables
 of them), we have, for example, three toric degenerations,
 whose toric varieties are given by 
\begin{enumerate}
\item $xy = zw,  zw = uv$,
\item $x^2=zw, zw=uv$,
\item $x^2=yz, w^2=uv$,
\end{enumerate}
 here $u, v, w, x, y, z$ are homogeneous coordinates of $\Bbb P^5$. 
The moment polytopes of them are 
 octahedron, quadrangular pyramid and tetrahedron respectively.
Jeffrey and Weitsman \cite{JW}
 explicitly described the Goldman system in this case and showed that the 
 moment polytope will be either a quadrangular pyramid or a tetrahedron,
 depending on the pants decomposition of the genus two curve.
However, the toric varieties associated to these polytopes do not
 satisfy Assumption 2, so it is not very good from the enumerative point of view.

Here we investigate the case of octahedron.
The torus action is given by
\[
 [x : y : z : w : u : v]
  \mapsto
   [\alpha x : \beta y : \gamma z : \alpha \beta \gamma^{-1} w
     : \alpha \beta u : v],
\alpha, \beta, \gamma\in\Bbb C^*.\]
The moment polytope is the convex hull of
\[
 \{ (\lambda, 0, 0), (0, \lambda, 0), (0, 0, \lambda),
    (\lambda, \lambda, - \lambda), (\lambda, \lambda, 0), (0, 0, 0) \},
  \lambda\geq 0,\] 
 and the defining inequalities are given by
\[
\begin{array}{ll}
 \ell_1(u) = \langle (0, 1, 1), u \rangle \ge 0,&
 \ell_2(u) = \langle (-1, 0, 0), u \rangle + \lambda \ge 0,\\
 \ell_3(u) = \langle (0, -1, 0), u \rangle + \lambda \ge 0,&
 \ell_4(u) = \langle (1, 0, 1), u \rangle \ge 0,\\
\ell_5(u) = \langle (0, 1, 0), u \rangle \ge 0, &
 \ell_6(u) = \langle (-1, 0, -1), u \rangle + \lambda \ge 0,\\
 \ell_7(u) = \langle (0, -1, -1), u \rangle + \lambda \ge 0, &
 \ell_8(u) = \langle (1, 0, 0), u \rangle \ge 0.
 \end{array}\]
It is easy to see that $X_0$ has a small resolution.
In fact, all the singularities of $X_0$
 are locally isomorphic to the singularity of the degeneration of $F_3$.
As before, 
 by the Fredholm regularity of torically transversal rational curve in toric varieties
 (\cite{CO}), we have the following.
\begin{prop}\label{curvelift}
Any torically transverse rational stable map
 in $X_0$ can be deformed to a rational stable map
 in $X$.  \qed
\end{prop} 
By Proposition \ref{homgr}, it is easy to see that in this case there is a
 natural isomorphism 
\[
H_2(X, \Bbb Z)\cong H_2(X_0, \Bbb Z).
\] 
One sees that there are four families of $\Bbb P^1$s on $X_0$
 corresponding to the pairs of parallel facets of the octahedron.
These curves have two dimensional freedom (given by the parallel transport)
 to move in $X_0$.
These are all homologous, and generates $H_2(X_0, \Bbb Z)$. 
By the isomorphism above, 
 any of the lifts of them due to Proposition \ref{curvelift}
 generates $H_2(X, \Bbb Z)$.
 
Now since the singular points of $X_0$ are locally isomorphic to  
 that of the degeneration of 
 $F_3$, the inverse images of them under the gradient flow $\phi_{grH, t}$
 are three dimensional spheres.
On the other hand, $X$ is ruled by any one of the two dimensional families of
 $\Bbb P^1$ which are the lifts of the rulings on $X_0$.
General members of these rulings
 on $X_0$ do not  intersect with the singular points,
 so general members of the corresponding rulings on $X$
 do not intersect with the three dimensional spheres above. 

Let $p_i$, $i = 1, \dots, a$ be general points on $X_0$ and
 $l_j$, $i = 1, \dots, b$ be general lines in one of the rulings of $X_0$,
 $a, b\in \Bbb Z_{\geq 0}$.
Extend them to a holomorphic family $p_{i, t}, l_{j, t}$ on $\mathfrak X\to \Bbb C$
 around the origin of $\Bbb C$. 
The following is a consequence of Propositions
 \ref{lift} and \ref{prop:converse}.

\begin{prop}\label{curvedeg}
There is a natural one to one correspondence between 
 the families of rational stable maps 
 \[\varphi: C_t\to X_t\] with 
 \[c_1(X_t)(\varphi(C)) = 
 2a+b\] intersecting $p_{i, t}, l_{j, t}$ for $t\neq 0$,
 and the rational stable maps 
 \[\psi: C'\to X_0\] with
 \[c_1(X_0)(\psi(C')) = 2a+b\] intersecting
 $p_i, l_j$.
The domains $C_t, C'$ are in fact $\Bbb P^1$, and $\psi$ is torically transverse.\qed
\end{prop}

\begin{rem}
In the case of degeneration to the quadrangular pyramid or to the tetrahedron,
 the rulings of $X_t$ are broken and it seems difficult to calculate the
 invariants of $X_t$ directly from the tropical calculation on $X_0$.
Also, as we mentioned above, $X_0$ does not satisfy Assumption 2,
 and the behavior of holomorphic curves
  may change when we move from $X_t$ to $X_0$.
This is the reason that Goldman systems may not be good from the view point of
 toric degenerations of integrable systems
 (or from the view point of enumerative geometry, either).
 
Concerning this point, it is an interesting problem to find a natural construction of a
 structure of an integrable
 system on $X$ which torically degenerate to $X_0$
 (note that we can pull-back the toric integrable system structure on $X_0$ to
 $X$ by $\phi_{grH}^{-1}$. 
 However, its geometric meaning is not clear, compared to Goldman systems).
\end{rem}
By Propositions \ref{curvelift} and \ref {curvedeg}, 
 we can compute the Gromov-Witten invariants of $X$
 by counting curves in $X_0$, 
 which in turn can be calculated by tropical method.
In this case, $H_2(X, \Bbb Z)$ is 
 free of rank one and so we can parameterize it
 by the set of integers.
Let $l$ be the generator of $H_2(X, \Bbb Z)$
 represented by a rational stable map.
Let $h$ be the generator of $H_4(X, \Bbb Z)$ such that 
 the intersection pairing satisfies
 \[l\cdot h = 1.\]

In this case, the set $\Delta_{a\cdot l}$ of coarse-degrees, $a\in \Bbb N$ 
 is described as follows.
\begin{lem}
A coarse-degree $\widetilde{\Delta}$ belongs to $\Delta_{a\cdot l}$
 if and only if the following conditions are satisfied:
 
Representing a coarse-degree as an integer valued function on the set of 
 primitive vectors in $N = \Bbb Z^3$, 
\[\begin{array}{l}
\widetilde{\Delta}(0, 1, 1) =  \widetilde{\Delta} (0, -1, -1) = a_1,\\
\widetilde{\Delta}(1, 0, 0) = \widetilde{\Delta} (-1, 0, 0) = a_2,\\
\widetilde{\Delta}(0, 1, 0) =  \widetilde{\Delta} (0, -1, 0) = a_3,\\
\widetilde{\Delta}(1, 0, 1) =  \widetilde{\Delta} (-1, 0, -1) = a_4,
\end{array}
\]
 where $a_1+a_2+a_3+a_4 = a$.
 The others are zero.\qed
\end{lem}
Then we have the following.
\begin{prop}
When \[2m = 2a+b,\] 
then
\[I_{0, k, m\cdot l}((pt^{\vee})^a, (l^{\vee})^b, (h^{\vee})^c) = 
 m^cI_{0, k', m\cdot l}((pt^{\vee})^a, (l^{\vee})^b),\]
 where $k = a+b+c$, $k' = a+b$.
The number
 \[I_{0, k', m\cdot l}((pt^{\vee})^a, (l^{\vee})^b)\]
 is the same as the number of $k'$-marked
 rational tropical curves of degree $\Delta$ incident to
  generic $a$ points and generic $b$ lines parallel to 
  either $(1, 0, 0), (0, 1, 0), (0, 1, 1)$ or $(1, 0, 1)$.
The degree $\Delta$ must satisfy the condition that the associated 
 coarse-degree $\Delta$ is contained in $\Delta_{m\cdot l}$.
When $2m\neq 2a + b$, then 
 \[I_{0, k, m\cdot l}((pt^{\vee})^a, (l^{\vee})^b, (h^{\vee})^c) = 0.\]\qed 
\end{prop}
\begin{rem}\label{rem:inv}
It might be interesting that the counting number of the tropical curves in the proposition 
 is the same however we choose the number of lines parallel to 
 $(1, 0, 0), (0, 1, 0), (0, 1, 1)$ or $(1, 0, 1)$, provided the total number is $b$.
This is not a purely tropical consequence, but follows from the 
 homological invariance of the Gromov-Witten invariants.
\end{rem}
Using our method, we can also calculate the Gromov-Witten invariants with 
 odd degree arguments.
First, as we remarked above, the generators $\sigma_i$, $i = 1, 2, 3, 4$
  of $H_3(X)$ are collapsed to  singular points of $X_0$
 by the map $\phi_{grH}$.
If some $I_{0, l, \beta}\neq 0$ with $\sigma_i$ in the argument, 
 then there is a family of rational stable maps
  $\psi_t: C_t\to X_t$, $t\in (0, 1]$
 satisfying the incidence conditions.
Assume that there are other arguments other than $\sigma_i$s.
By Gromov's compactness theorem, 
 we can assume that $\psi_t$ converges to a rational stale map in
 $X_0$.
 
By the above remark about the classes $\sigma_i$, 
 the limit curve must intersect the singular point of $X_0$.
On the other hand, it satisfies the incidence conditions induced from 
 the classes other than $\sigma_i$.
However, as in the proof of Lemma \ref{transverse}, such a curve belongs to a lower
 dimensional subvariety of the moduli space, 
 and so it cannot satisfy these incidence conditions.
The case where all the arguments are from $\sigma_i$s can be dealt with
 by similar dimension counting argument.
So we have the following.
\begin{prop}
The number $I_{0, l, \beta}(\alpha_1, \dots, \alpha_k) = 0$
 if one of $\alpha_i$ is a cohomology class of odd degree.\qed
\end{prop}

%%%%%%%%%%%%%%%%%%%%%%%%%%%%%%%%%%%%%%%%%%%%%%%%%%%%%%%%%%%%%%%%%%%%%%
 
%%%%%%%%%%%%%%%%%%%%%%%%%%%%%%%%%%%%%%%%%%%%%%%%%%%%%%%%%%%%%%%%%%%%%%


\begin{thebibliography}{99} 
\bibitem{CO}{\sc C. Cho and Y. Oh},
     {\it Floer cohomology and disc instantons of {L}agrangian torus
              fibers in {F}ano toric manifolds},
   {Asian J. Math.},
    {10},
       {(2006)},
   {no. 4},
    {773--814}.
\bibitem{CK}{\sc  D. Cox and S. Katz,}
 {\it Mirror symmetry and algebraic geometry.}
 Mathematical Surveys and Monographs, 68.
 American Mathematical Society, Providence, RI, 1999. 
\bibitem{G}{\sc W. Goldman,}
{\it Invariant functions on Lie groups and Hamiltonian
flows of surface group representations.}
 Invent. Math., 85, 1986, no. 2, 263--302.
\bibitem{JW}{\sc L. Jeffrey and J. Weitsman,}
{\it Bohr-Sommerfeld orbits in the
moduli space of flat connections and the Verlinde dimension formula.}
Comm. Math. Phys., 150, (1992), no. 3, 593--630.
%
\bibitem{KM}{\sc M. Kogan, and E. Miller},
     {\it Toric degeneration of {S}chubert varieties and
              {G}elfand-{T}setlin polytopes},
   {Adv. Math.},
      {193},
      {2005},
    {no. 1},
     {1--17},
\bibitem{M}{\sc G. Mikhalkin,}
{\it Enumerative tropical algebraic geometry in $\Bbb R\sp 2$.}
  J. Amer. Math. Soc.  18  (2005),  no. 2, 313--377
\bibitem{NR}{\sc M. Narasimhan and S. Ramanan,}
{\it Moduli of vector bundles on a compact
Riemann surface.} Ann. of Math. 89, (1969), no. 2, 14?51.
\bibitem{N}{\sc P. Newstead,}
{\it Stable bundles of rank 2 and odd degree over a curve of
genus 2.} Topology, 7 (1968), 205?215.
\bibitem{N1}{\sc T. Nishinou,}
{\it Disc counting on toric varieties via tropical curves},
 preprint.
\bibitem{N2}{\sc T. Nishinou,}
{\it Correspondence theorems for tropical curves}, preprint.
\bibitem{NNU}{\sc T. Nishinou et al.,}
{\it Toric degenerations of Gelfand-Cetlin systems and potential
functions}, Adv. Math. (2009), doi:10.1016/j.aim.2009.12.012.
\bibitem{NS}{\sc T. Nishinou and B. Siebert,}
{\it Toric degenerations of toric varieties and tropical curves.}
 Duke Math. J. 135 (2006), no. 1, 1--51.
\bibitem{R}{\sc W. Ruan},
      {\it Lagrangian torus fibration of quintic {C}alabi-{Y}au
              hypersurfaces. {II}. {T}echnical results on gradient flow
              construction},
    {J. Symplectic Geom.},
   {1},
    {2002},
 {no. 3},
     {435--521}.
\end{thebibliography}
\end{document}